\documentclass[12pt]{article}
\usepackage[latin1]{inputenc}
\usepackage{amssymb}
\usepackage{amsmath}
\usepackage{amsthm}

\input pictex.tex
\newcommand{\esp}{\hspace{0.1cm}}

\newcommand{\R}{\mathbb{R}}
\newcommand{\Q}{\mathbb Q}
\newcommand{\N}{\mathbb N}
\newcommand{\Z}{\mathbb Z}

\theoremstyle{definition}

\newtheorem{thm}{Theorem}[section]

\newtheorem{prop}[thm]{Proposition}

\newtheorem{lem}[thm]{Lemma}

\newtheorem{rem}[thm]{Remark}

\newtheorem{cor}[thm]{Corollary}

\newcommand{\vs}{\vspace{0.4cm}}

\newcommand{\vsp}{\vspace{0.1cm}}

\date{}
\author{}

\begin{document}

\title{On groups with finitely many Conradian orderings}
\author{Crist\'obal Rivas}
\maketitle
\begin{abstract} We study the space of left-orderings
on groups with (only) finitely many Conradian orderings. We show
that, within this class of groups, having an isolated left-ordering
is equivalent to having finitely many left-orderings.
\end{abstract}

\vs


\vsp


\vs

\noindent \textbf{\large Introduction}

\vs

A (total) left-ordering $\preceq$ on a group $G$ is said to be
\textit{isolated} if there is a finite family
$\{g_1,\ldots,g_n\}\subset G$ such that $\preceq$ is the only
left-ordering on $G$ with the property that $g_i\succ id$, for
$1\leq i \leq n$. This criteria may be used to define a topology on
$\mathcal{LO}(G)$, the set of all left-orderings on $G$. It was proved by Sikora in
\cite{sikora} that with this topology, $\mathcal{LO}(G)$ is a
totally disconnected, Hausdorff and compact topological space.
Moreover, when $G$ is countable, this topology is metrizable. See \S \ref{space of orderings} for further details.

\vsp

Knowing whether a given group has an isolated left-ordering has been
a question of major interest in the recent development of the theory
of orderable groups. A big progress was made by Tararin who
classified left-orderable groups that admit only finitely many
left-orderings (a \textit{Tararin group}, for short), see Theorem
\ref{teo T} or \cite[\S 5.2]{kopytov}.

\vsp

Albeit Tararin's description has shown to be very useful, the
comprehension of groups admitting isolated left-orderings is far
from being reached. Some progress in this direction was done in
\cite{dd} and \cite{navas-hecke}. In \cite{dd}, Dubrovina and
Dubrovin show that braid groups have isolated left-orderings,
whereas in \cite{navas-hecke}, Navas describes a family of
two-generated groups (which contains the three strands braid group
$B_3$) having infinitely many left-orderings together with isolated
left-orderings. For a nice survey about orderings on braid groups,
see \cite{braids}.

\vsp

It follows from Tararin's description that every Tararin group is
solvable. On the other hand, neither braid groups nor the groups
described in $\cite{navas-hecke}$ are solvable. Moreover, in
\cite{navas} it is shown that the only nilpotent groups having
isolated left-orderings are the torsion-free, rank-one Abelian
groups\footnote{Recall that a torsion-free Abelian group $\Gamma$
has \textit{rank} $n$ if $n$ is the least integer for which $\Gamma$ embeds into $\Q^n$.}. Thus, it is natural to pose the

\vsp\vsp

\noindent \textbf{Main Question:} Is it true that, in the class of
left-orderable solvable groups, having an isolated left-ordering is
equivalent to having only finitely many left-orderings?

\vsp\vsp

In this work we give a partial (affirmative) answer to this
question. Recall that a left-ordering $\preceq$ on a group $G$ is
\textit{Conradian} (or a $C$-ordering) if $f\succ id$ and $g\succ
id$ imply $fg^2\succ g$, see \cite{conrad,leslie, navas}.

\vs

\noindent \textbf{Main Theorem:} \textit{Let $G$ be a group
admitting only finitely many $C$-orderings. Then $G$ either admits
only finitely many left-orderings (so $G$ is a Tararin group) or has
no isolated left-orderings.}

\vs

We note that the relation between left-orderings and Conradian orderings is much deeper than just the one described in the Main Theorem. For instance, in \cite[\S 4]{navas} it is proved that no Conradian ordering is isolated in a group with infinitely many left-orderings, and also a criterion is given for a left-ordering to be isolated in terms of the so-called {\em Conradian soul} of an ordering. Nevertheless, we will not make use of those facts in this work.

\vsp

To prove the  Main Theorem we will make use of the algebraic description of
groups admitting (only) finitely many $C$-orderings, here Theorem
\ref{teo B}, which was obtained in \cite{rivas}. As shown in Theorem
\ref{teo B}, groups with finitely many Conradian orderings admits a unique rational series (see definition below), and our proof proceeds by induction on the length of this series. In \S 2, we explore the (initial) case of groups with rational series of length two. In
this case, we give an explicit description of
$\mathcal{LO}(G)$. In \S 3.1 we obtain some
technical results concerning the action of inner automorphisms of a
group $G$ with a finite number of Conradian orderings. As a
consequence, we show that the maximal convex subgroup of $G$ (with
respect to a $C$-ordering) is a group that fits into the
classification made by Tararin. Finally, in \S 3.2, we prove the
general case, while \S 3.3 is devoted to the description of an
illustrative example.

%%%%%%%%%%%%%%%%%%%%%%%%%%%%%%%%%%%%%%%%%%%%%%%%%%%%%

\section{Preliminaries}

%%%%%%%%%%%%%%%%%%%%%%%%%%%%%%%%%%%%%%%%%%%%%%%%%%%%%%%

We begin this section recalling the foundational result
\cite[Theorem 4.1]{conrad}. Recall that in a left-ordered group $G$,
$G_g$ (resp. $G^g$) denotes the maximal (resp. minimal) convex
subgroup which does not contain (resp. contains) $g\in G$. (A subset
$S$ of a left-ordered group $\Gamma$ is said to be convex if and
only if for every $\gamma\in \Gamma$ such that $s_1\preceq \gamma
\preceq s_2$, for some $s_1,\; s_2$ in $S$, we have that $\gamma\in
S$.)

\begin{thm}[Conrad] \label{teo C}{\em An ordering $\preceq$ on a group $G$ is Conradian
if and only if for every $g\in G$, $g\not=id$, we have that $G_g$ is
normal in $G^g$, and there exists a unique up to multiplication by
a positive real number, non-decreasing group homomorphism  \esp
$\tau_{\preceq}^{g} \!: G^g \rightarrow \mathbb{R}$ \esp whose
kernel coincides with $G_g$.}
\end{thm}

The Conrad Theorem implies that any $C$-orderable group is locally
indicable\footnote{A group $G$ is locally indicable if for any
finitely generated subgroup $H$ there is a non-trivial group
homomorphism from $H$ to the real numbers under addition.}, and a
remarkable result from \cite{brodskii} shows that the class of
$C$-orderable groups coincides with the class of locally indicable
groups, see also \cite{navas}. Thus, all torsion-free, one-relator
groups are $C$-orderable \cite{brodskii, howie}.

\vsp

In \cite{rivas} a structure theorem was given for groups admitting
only finitely many Conradian orderings. For the statement, recall
that a series
$$\{ id \} = G_0 \lhd G_{1} \lhd \ldots \lhd G_{n-1} \lhd G_n = G$$
is said to be {\em rational} if it is subnormal ({\em i.e.,} each
$G_{i}$ is normal in $G_{i+1}$) and each quotient $G_{i+1} / G_{i}$
is torsion-free rank-one Abelian. We say that the rational series is
{\em normal} if, in addition, $G_i\lhd G$ for all $1\leq i \leq n$.

\begin{thm} \label{teo B} {\em Let $G$ be a $C$-orderable group. If $G$
admits only finitely many $C$-orderings, then $G$ admits a unique
(hence normal) rational series. In this series, no quotient
$G_{i+2}/G_{i}$ is Abelian. Conversely, if $G$ is a group admitting
a  normal rational series
$$\{ id \} = G_0 \lhd G_{1} \lhd \ldots \lhd G_{n-1} \lhd G_n = G$$
so that no quotient $G_{i+2} / G_{i}$ is Abelian, then the number of
$C$-orderings on $G$ equals $2^n$.}
\end{thm}

One of the crucial steps in proving Theorem \ref{teo B} consist in using the
Conrad Theorem to show that in any C-ordering of $G$ -a
group with only finitely many Conradian orderings- and any $g\in G$,
we have that $G_g=G_i$ and $G^g=G_{i+1}$ for some $0\leq i\leq n-1$.
In particular, in a group with only finitely many Conradian
orderings, the convex series given by a C-ordering coincides
with the rational series of $G$.

\vsp

A sub-class of the class of groups admitting only finitely many
Conradian orderings is the class of groups admitting only finitely
many left-orderings. This latter class was described by Tararin,
\cite[\S 5.2]{kopytov}. Since we will make use of this description,
we quote Tararin's theorem below. For the statement, recall that a left-ordering $\preceq$ on a group G is said to be {\em bi-invariant} (or {\em bi-ordering}, for short) if $g\succ id$ implies $hgh^{-1}\succ id$ for all $h\in G$. Clearly, every bi-ordering is Conradian.

\begin{thm}[Tararin] \label{teo T}{\em Let $G$ be a left-orderable group. If $G$
admits only finitely many left-orderings, then $G$ admits a unique
(hence normal) rational series. In this series, no quotient
$G_{i+2}/G_{i}$ is bi-orderable. Conversely, if $G$ is a group
admitting a normal rational series
$$\{ id \} = G_0 \lhd G_{1} \lhd \ldots \lhd G_{n-1} \lhd G_n = G$$
so that no quotient $G_{i+2} / G_{i}$ is bi-orderable, then the
number of left-orderings on $G$ equals $2^n$.}
\end{thm}

Note that the statement of Tararin's theorem is the same as the
statement of Theorem \ref{teo B} though changing `$C$-orderings' by
`left-orderings', and the condition `$G_{i+2} / G_{i}$ non-Abelian'
by `$G_{i+2} / G_{i}$ non-bi-orderable'.

\vsp

%%%%%%%%%%%%%%%%%%%%%%%%%%%%%%%%%%%%%%%%%%%%%

\subsection{The space of left-orderings of a group}
\label{space of orderings}

%%%%%%%%%%%%%%%%%%%%%%%%%%%%%%%%%%%%%%%%%%%%%%%%%%

Recall that given a left-ordering $\preceq$ on a group $G$, we say
that $f\in G$ is {\em $\preceq$-positive} or simply \textit{positive}
(resp. {\em $\preceq$-negative} or \textit{negative}) if $f\succ id$
(resp. $f\prec id$). We denote $P_\preceq$ the set of
$\preceq$-positive elements in $G$. Clearly, $P_\preceq$ satisfies the
following properties:

\vsp

\noindent $(i)$ $P_\preceq  P_\preceq \subseteq P_\preceq \;$, that
is, $P_\preceq$ is a semigroup;

\vsp

\noindent $(ii)$ $G=P_\preceq \sqcup P^{-1}_\preceq \sqcup \{id\}$,
where the union is disjoint, and $P^{-1}_\preceq=\{g^{-1}\in G\mid
g\in P_\preceq\}=\{g\in G\mid g\prec id\}$.

\vsp

Moreover, given any subset $P\subseteq G$ satisfying the conditions
$(i)$ and $(ii)$ above, we can define a left-ordering $\preceq_P$ by
$f\prec_P g$ if and only if $f^{-1}g\in P$. Therefore, describing a left-ordering is equivalent to describing its set of positive elements. We usually identify $\preceq$ with $P_{\preceq}$.

\vsp

Given a left-orderable group $G$ (of arbitrary cardinality), we
denote the set of all left-orderings on $G$ by $\mathcal{LO} (G)$.
This set has a natural topology first introduced by Sikora for the case of countable groups \cite{sikora}. This topology can be defined by identifying $P\in
\mathcal{LO}(G)$ with its characteristic function $\chi_P \in
\{0,1\}^G $. In this way, we can view $\mathcal{LO}(G)$ embedded in
$\{0,1\}^G$. This latter space, with the product topology, is a
Hausdorff, totally disconnected, and compact space. It is not hard
to see that (the image of) $\mathcal{LO}(G)$ is closed inside, and
hence compact as well (see \cite{navas,sikora} for details).

\vsp

A basis of neighborhoods of $\,\preceq \,$ in $\, \mathcal{LO} (G)$ is the family of the sets $\, V_{f_1,\ldots,f_k} \,$ of all left-orderings $\, \preceq' \,$ on $G$ such that all the $\,f_i\,$ are
$\,\preceq'$-positive, where $\{f_1,\ldots,f_k\}$ runs over all
finite subsets of $\,\preceq$-positive elements of $G$. Hence, a left-ordering of $G$ is {\em isolated} (in the sense of the introduction) if an only if it is an isolated point of $\mathcal{LO}(G)$. The (perhaps
empty) subspaces $\, \mathcal{BO}(G) \,$ and $\, \mathcal{CO}(G) \,$
of bi-orderings and $C$-orderings on $G$ respectively, are closed
inside $\mathcal{LO}(G)$, hence compact; see \cite{navas}.

\vsp

If $G$ is countable, then this topology is metrizable: given an
exhaustion $G_0 \subset G_1 \subset \ldots$ of $G$ by finite sets,
for different $\, \preceq \,$ and $\,\preceq'\,, \,$ we may define
$\,dist(\preceq,\preceq') = 1 / 2^n$, where $n$ is the first integer
such that $\, \preceq \,$ and $\, \preceq'\,$ do not coincide on
$G_n$. If $G$ is finitely generated, we may take $G_n$ as the ball
of radius $n$ with respect to a fixed finite system of generators.

%%%%%%%%%%%%%%%%%%%%%%%%%%%%%%%%%%%

\subsection{A basic construction for producing new left-orderings}

%%%%%%%%%%%%%%%%%%%%%%%%%%%%%%%%%%%

In this section we describe some basic constructions for creating
new left-orderings starting with a given one. The main
idea is to exploit the flexibility given by the convex subgroups.

\vsp

Let $\preceq$ be a left-ordering on a group $G$. If $C$ is a proper convex subgroup of $G$, then $\preceq$ induces a total order $\preceq^C$ on the set of
left-cosets of $C$ by
\begin{equation} \label{def convex
ord}g_1C\prec^{C} g_2C \Leftrightarrow g_1c_1\prec g_2c_2\;\; \text{for
all $c_1, c_2$ in C}.\end{equation} More importantly, this order is
preserved by the left action of $G$; see for instance \cite[\S
2]{kopytov}. In particular, if $C$ is a normal subgroup, then $\preceq^C$ becomes a left-ordering of the group $G/C$.

\vsp

As the reader can easily check, the left-ordering $\preceq$ can
be recovered from the left-ordering $\preceq^{C}$ and the
left-ordering $\preceq_C$, defined as the restriction of $\preceq$
to $C$, by the following equation:
$$ g\succ id \Leftrightarrow \left\{ \begin{array}{l } gC\succ^C C  \text{ or } \\ gC=C \;\;\text{ and } g\succ_C id . \end{array}
\right. $$

This easily implies

\begin{lem} {\em Let $\preceq$ be a left-ordering on a group $G$, and suppose there is a non-trivial convex subgroup $C$. Then there is a continuous injection $$\varphi:\mathcal{LO}(C)\to \mathcal{LO}(G)$$ such that $\preceq$ belongs to the image of $\varphi$.

\vsp

Moreover, if in addition $C$ is normal, then we have a continuous injection $$\varphi:\mathcal{LO}(G/C)\times \mathcal{LO}(C)\to\mathcal{LO}(G)$$ such that $\preceq$ belongs to the image of $\varphi$.}

\end{lem}

\begin{cor} \label{extension} {\em If for a left-ordering $\preceq$ on a group $G$ there is a convex subgroup $C$ such that either $C$ has no isolated left-orderings or such that $C$ is normal and $G/C$ has no isolated left-orderings, then $\preceq$ is non-isolated.}
\end{cor}

%%%%%%%%%%%%%%%%%%%%%%%%%%%%%%%%%%%%%%%%%%%%%%%%%%%%%

\section{On groups with a rational series of length two}

%%%%%%%%%%%%%%%%%%%%%%%%%%%%%%%%%%%%%%%%%%%%%%%%%%%%%%%

\hspace{0.4 cm} Throughout this section, $G$ will denote a
left-orderable, non-Abelian group with a rational series of length
$2$:
$$\{id\}=G_0\lhd G_1\lhd G_2=G.$$

\vsp

If the group $G$ is not bi-orderable, then $G$ has a normal rational
series of length 2 and the quotient $G_2/G_0=G$ is non-bi-orderable.
Thus $G$ fits into the classification made by Tararin, so it has
only finitely many left-orderings.

\vsp

For the rest of this section we will assume that $G$ is not a
Tararin group, so $G$ is bi-orderable. We have

\begin{lem} \label{lema 1}\textit{The group $G$ satisfies that $G/G_1\simeq\Z$.}
\end{lem}

\noindent \textit{Proof:} Consider the action by conjugation
$\alpha: G/G_1\to Aut(G_1)$ given by $\alpha(gG_1)(h)=ghg^{-1}$.
Since $G$ is non-Abelian, we have that this action is non-trivial,
{\em i.e.} $Ker(\alpha)\not=G/G_1$. Moreover, $Ker(\alpha)=\{id\}$,
since in the other case, as $G/G_1$ is rank-one Abelian, we would
have that $(G/G_1)/Ker(\alpha)$ is a torsion group. But the only
non-trivial, finite order automorphism of $G_1$ is the inversion,
which implies that $G$ is non-bi-orderable, thus a Tararin group.

\vsp

The following claim is elementary and we leave its proof to the
reader.

\vsp

\noindent \underbar{Claim:} If $\Gamma$ is a torsion-free, rank-one
Abelian group such that $\Gamma \not\simeq  \Z$, then for any $g\in
\Gamma$, there is an integer $n>1$ and $g_n\in \Gamma$ such that
$g_n^n=g$.

\vsp

Now take any $b\in G\setminus G_1$ so that $\alpha (bG_1)$ is a
non-trivial automorphism of $G_1$. Since $G_1$ is rank-one Abelian,
for some positive $r=p/q\in \Q$, $r\not=1$, we must have that
$bab^{-1}=a^r$ for all $a\in G_1$. Suppose that $G/G_1\not\simeq
\Z$. By the previous claim, we have a sequence of increasing
integers $(n_1,n_2\ldots)$ and a sequence $(g_1,g_2,\ldots)$ of
elements in $G/G_1$ such that $g_i^{n_i}=bG_1$. In particular we have that $g_iag_i^{-1}=a^{r_i}$, where $r_i$ is a rational such that $r_i^{n_i}=r$. In other words, given $r$, we have found among the rational numbers, an infinite collection of $r_i$ solving the
equation $x^{n_i}-r=0$, but, by the Rational Roots Theorem or
Rational Roots Test \cite[Proposition 5.1]{morandi}, this can not
happen. This finishes the proof of Lemma \ref{lema 1}. $\hfill \square$

\vs

\begin{lem}\label{lema afin} {\em The group $G$ embeds in $Af_+(\R)$, the group of (orientation preserving)
affine homeomorphism of the real line.}
\end{lem}

\noindent \textit{Proof:} We first embed $G_1$. Fix $a\in G_1$,
$a\not=id$. Define $\varphi_a:G_1\to Af_+(\R)$ by declaring
$\varphi_a(a)(x)=x+1$, and if $a^\prime \in G_1$ is such that
$(a^\prime)^q=a^p$, we declare $\varphi_a(a^\prime)(x)=x+p/q$.
Showing that $\varphi_a$ is an injective homomorphism is routine.

\vsp

Now let $b\in G$ such that $\langle bG_1\rangle=G/G_1$. Let
$1\not=r\in \Q$ such that $ba^\prime b^{-1}=(a^\prime)^r$ for every
$a^\prime\in G_1$. Since $G$ is bi-orderable we have that $r>0$, and
changing $b$ by $b^{-1}$ if necessary, we may assume that $r>1$.
Then, given $w\in G$, there is a unique $n\in \Z$ and a unique
$\overline{w}\in G_1$ such that $w=b^n\overline{w}$.

\vsp

Define $\varphi_{b,a}:G\to Af_+(\R)$ by
$\varphi_{b,a}(b^n\overline{w}):=H^{(n)}_r\circ
\varphi_a(\overline{w}),$ where $H_r(x):=rx$ , and $H_r^{(n)}$ is
the $n$-th composition of $H_r$ (by convention $H_r^{(0)}(x)=x$). We
claim that $\varphi_{b,a}$ is an injective homomorphism.

\vsp

Indeed, let $w_1,w_2\in G$, $w_1=b^{n_1}\overline{w}_1$,
$w_2=b^{n_2}\overline{w}_2$. Let $r_1\in \Q$ be such that
$\varphi_a(\overline{w}_1)(x)=x+r_1$. Then $ H_r^{(n)}\circ
\varphi_a(b^{-n}\overline{w}_1b^{n})(x)=H_r^{(n)} \circ
\varphi_a(\overline{w}_1^{(1/r)^{n}})(x)=r^n(x+
(1/r)^{n}r_1)=\varphi_a(\overline{w}_1)\circ H_r^{(n)}(x),$ for all
$n\in \Z$. Thus
\begin{eqnarray*}
\varphi_{b,a}(w_1w_2)
 &=& \varphi_{b,a}(b^{n_1}b^{n_2} \; b^{-n_2}\overline{w}_1 b^{n_2} \overline{w}_2)=H_r^{(n_1)} \circ H_r^{(n_2)}\circ \varphi_a(b^{-n_2}\overline{w}_1 b^{n_2})
 \circ\varphi_a( \overline{w}_2)\\
 &=& H_r^{(n_1)} \circ \varphi_a(\overline{w}_1)\circ H_r^{(n_2)}
 \circ\varphi_a( \overline{w}_2)=\varphi_{b,a}(w_1)\circ \varphi_{b,a}(w_2).
\end{eqnarray*}
So $\varphi_{b,a}$ is a homomorphism. To see that it is injective, suppose
that $\varphi_{b,a}(w_1)(x)=\varphi_{b,a}(b^{n_1}\overline{w}_1)(x) =
r^n\,x+r^nr_1=x$ for all $x\in \R$. Then $n=0$ and $r_1=0$, showing
that $w_1=id$. This finishes the proof of Lemma \ref{lema afin}. $\hfill\square$

\vs

Once the embedding $\varphi:=\varphi_{b,a}:G\to Af_+(\R)$ is fixed,
we can associate to each irrational number $\varepsilon$ an {\em
induced left-ordering} $\preceq_{\varepsilon}$ on $G$ whose set of
positive elements is defined by $\{g \!\in\! G \mid \esp \varphi
(g)(\varepsilon)
> \varepsilon \}$. When $\varepsilon$ is rational, the
preceding set defines only a partial ordering. However, in this case
the stabilizer of the point $\varepsilon$ is isomorphic to
$\mathbb{Z}$, and hence this partial ordering may be completed to
two total left-orderings $\preceq_{\varepsilon}^+$ and
$\preceq_{\varepsilon}^{-}$. These orderings were introduced by
Smirnov in \cite{smirnov}. Once the representation $\varphi$ is
fixed, we call these orderings, together with its corresponding
reverse orderings, Smirnov-type orderings. (By definition the
reverse ordering of $\preceq$, denoted $\overline{\preceq}$,
satisfies $id\,\overline{\prec} \, g$ if and only if $id\prec
g^{-1}$.)

\vsp

Besides the Smirnov-type orderings on $G$, there are four Conradian (actually bi-invariant!)
orderings. Since $G_1$ is always convex in
a Conradian ordering, $b^na^s\in G$, $n\not=0$, is positive if and only if $b$ is positive.
Then it is not hard to check that the four Conradian orderings are the following:

\vsp\vsp

1) $\preceq_{C_1}$, defined by $id \prec_{C_1} b^n a^s$ ($n\in \Z$,
$s\in \Q$) if and only if $n\geq 1$, or $n=0$ and $s>0$.

\vsp

2) $\preceq_{C_2}$, defined by $id \prec_{C_2} b^na^s$ if and only
if $n\leq -1$, or $n=0$ and $s>0$.

\vsp

3) $\preceq_{C_3}=\overline{\preceq}_{C_1} $.

\vsp

4) $\preceq_{C_4}=\overline{\preceq}_{C_2}$.

\begin{prop} \label{C y S} {\em Let $U\subseteq \mathcal{LO}(G)$ be
the set consisting of the four Conradian orderings together with the
Smirnov-type orderings. Then any ordering in $U$ is non-isolated in
$U$.}
\end{prop}

\noindent {\em Proof:} We first show that the Conradian orderings
are non-isolated.

\vsp

We claim that $\preceq_\varepsilon \to \preceq_{C_1}$ when
$\varepsilon\to \infty$. For this it suffices to see that any
positive element in the $\preceq_{C_1}$ ordering becomes
$\preceq_\varepsilon$-positive for any $\varepsilon$ is large
enough.

\vsp

By definition of $\preceq_\varepsilon$ we have that
$id\prec_\varepsilon b^na^s$ if and only if
$r^n(\varepsilon+s)=\varphi(b^na^s)(\varepsilon)>\varepsilon$, where
$r>1$. Now, assume that $id\prec_{C_1}b^na^s$. If $n=0$, then $s>0$
and $\varepsilon+s>\varepsilon$. If $n\geq 1$, then
$r^n(\varepsilon+s)>\varepsilon$ for $\varepsilon
>\frac{-r^ns}{r^n-1}$. So the claim follows.

\vsp

For approximating the other three Conradian orderings, we first note
that, arguing just as before, we have $\preceq_\varepsilon \to
\preceq_{C_2}$ when $\varepsilon\to -\infty$. Finally, the other two
Conradian orderings $\overline{\preceq}_{C_1}$ and
$\overline{\preceq}_{C_2}$ are approximated by
$\overline{\preceq}_\varepsilon$ when $\varepsilon \to \infty$ and
$\varepsilon\to -\infty $ respectively.

\vsp

Now let $\preceq_S$ be an Smirnov-type ordering and let
$\{g_1,\ldots, g_n\}$ be a set of $\preceq_S$-positive elements.

\vsp

Suppose first that $\preceq_S$ equals $\preceq_\varepsilon$, where
$\varepsilon$ has free orbit. Then we have that
$\varphi(g_i)(\varepsilon)>\varepsilon$ for all $1\leq i\leq n$.
Thus, if $\varepsilon^\prime$ is such that
$\varepsilon<\varepsilon^\prime< \min\{\varphi(g_i)(\varepsilon)\}$,
$1\leq i\leq n$, then we still have that
$\varphi(g_i)(\varepsilon^\prime)
>\varepsilon^\prime$, hence $g_i\succ_{\varepsilon^\prime} id$ for $1\leq i \leq n$.
To see that
$\preceq_{\varepsilon^\prime}\not=\preceq_{\varepsilon}$, first
notice that $\varphi(G_1)(x)$ is dense in $\R$ for all $x\in \R$. In
particular, taking $g\in G_1$ such that $\varepsilon< \varphi(g)(0)
<\varepsilon^\prime$, we have that
$\varphi(gb^ng^{-1})(\varepsilon)=\varphi(g)
(r^n\varphi(g)^{-1}(\varepsilon))=r^n\varphi(g)^{-1}(\varepsilon)+\varphi(g)(0)$.
Since $\varphi(g)^{-1}(\varepsilon)<0$ we have that for $n$ large
enough $gb^ng^{-1}\prec_\varepsilon id$. The same argument shows
that $gb^ng^{-1} \succ_{\varepsilon^\prime}id$. Therefore
$\preceq_{\varepsilon}$ and $\preceq_{\varepsilon^\prime}$ are
distinct.

\vsp

The remaining case is when $\preceq_S=\preceq_\varepsilon^\pm$. In
this case we can order the set $\{g_1,\ldots, g_n\}$ such that there
is $i_0$ with $\varphi(g_i)(\varepsilon)>\varepsilon$ for $1\leq
i\leq i_0$, and $\varphi(g_i)(\varepsilon)=\varepsilon$ for
$i_0+1\leq i \leq n$. That is $g_i\in Stab(\varepsilon)\simeq \Z$
for $i_0+1\leq i \leq n$. Let $\varepsilon^\prime>\varepsilon$.

\vsp

We claim that either $\varphi(g_i)(\varepsilon^\prime)>
\varepsilon^\prime$ for all $i_0+1\leq i \leq n$ or
$\varphi(g_i)(\varepsilon^\prime)< \varepsilon^\prime$ for all
$i_0+1\leq i \leq n$. Indeed, since $\varphi$ gives an affine
action, it can not be the case that a non-trivial element of $G$
fixes two points. So we have that $\varphi(g_i)(\varepsilon^\prime)
\not=\varepsilon^\prime$ for each $i_0+1\leq i \leq n$. Now, suppose
for a contradiction that there are $g_{i_0}, \; g_{i_1}\in
Stab(\varepsilon)$ with $g_{i_0}(\varepsilon^\prime)
<\varepsilon^\prime$ and $g_{i_1}(\varepsilon^\prime)
>\varepsilon^\prime$. Let $n,m\in \N$ be such that $g_{i_0}^n
=g_{i_1}^m$. Then $\varepsilon^\prime
<\varphi(g_{i_1})^m(\varepsilon^{\prime}) = \varphi(g_{i_0})^n
(\varepsilon^{\prime})<\varepsilon^\prime$. A contradiction. So the
claim follows.

\vsp

Now assume that $\varphi(g_i)(\varepsilon^\prime)>
\varepsilon^\prime,$ for all $i_0+1\leq i\leq n$. If, in addition
$\varepsilon<\varepsilon^\prime< \min\{\varphi(g_i)(\varepsilon)\}$
with $1\leq i\leq i_0$, then $g_i\succ_{\varepsilon^\prime} id$ for
$1\leq i \leq n$, showing that $\preceq_S$ is non-isolated. In the
case where $\varphi(g_i)(\varepsilon^\prime)< \varepsilon^\prime$
for all $i_0+1\leq i\leq n$, we let $\tilde{\varepsilon}$ such that
$\max \{\varphi(g_i)^{-1}(\varepsilon)\}< \tilde{\varepsilon}
<\varepsilon$ for $1\leq i \leq i_0$. Then we have that
$g_i\succ_{\tilde{\varepsilon}} id$ for $1\leq i \leq n$. This shows
that, in any case, $\preceq_S=\preceq_\varepsilon^\pm$ is
non-isolated in $U$. $\hfill\square$

\vs

The following theorem shows that the space of left-orderings of $G$
is made up by the Smirnov-type orderings together with the Conradian
orderings. This generalizes \cite[Theorem 1.2]{rivas}.

\begin{thm}\label{laprop}\textit{Suppose $G$ is a non Abelian group with rational series of length 2.
If $G$ is bi-orderable, then its space of left-orderings has no isolated
points. Moreover, every non-Conradian ordering is equal to an
induced, Smirnov-type, ordering arising from an affine action of $G$
over $\R$ given by $\varphi$ above.}

\end{thm}

\vsp

To prove Theorem \ref{laprop}, we will use the ideas (and notation)
involved in the following well-known orderability criterion (see
\cite[Theorem 6.8]{ghys}, \cite[\S 2.2.3]{book}, or
\cite[Proposition 2.1]{navas} for further details).

\vsp

\begin{prop}\label{real din} {\em For a countable infinite group $\Gamma$, the following two properties are
equivalent:

\vsp

\noindent -- $\Gamma$ is left-orderable,

\vsp

\noindent -- $\Gamma$ acts faithfully on the real line by
orientation preserving homeomorphisms.}
\end{prop}

\noindent\textit{Sketch of proof: } The fact that a group of
orientation preserving homeomorphisms of the real line is
left-orderable is easy and may be found also in \cite[Theorem
3.4.1]{kopytov}. In what follows, we will not make use of this.

\vsp

For the converse, we construct what is called \textit{the dynamical
realization of a left-ordering}. Let $\preceq$ be a left-ordering on
$\Gamma$. Fix an enumeration $(g_i)_{i \geq 0}$ of $\Gamma$, and let
$t(g_0)=0$. We shall define an order-preserving map $t: \Gamma \to
\R$ by induction. Suppose that $t(g_0), t(g_1),\ldots,t(g_i)$ have
been already defined. Then if $g_{i+1}$ is greater (resp. smaller)
than all $g_0,\ldots, g_i$, we define $t(g_{i+1})=
max\{t(g_0),\ldots, t(g_i)\}+1$ (resp. $min\{t(g_0),\ldots,
t(g_i)\}-1$). If $g_{i+1}$ is neither greater nor smaller than all
$g_0,\ldots,g_i$, then there are $g_n,g_m\in\{g_0,\ldots , g_i \}$
such that $g_n\prec g_{i+1}\prec g_m$ and no $g_j$ is between
$g_n,g_m$ for $0\leq j\leq i$. Then we put
$t(g_{i+1})=(t(g_n)+t(g_m))/2$.

\vsp

Note that $\Gamma$ acts naturally on $t(\Gamma)$ by $g(t(g_i)) =
t(gg_i)$. It is not difficult to see that this action extends
continuously to the closure of $t(\Gamma)$. Finally, one can extend
the action to the whole real line by declaring the map $g$ to be
affine on each interval in the complement of $t(\Gamma)$.
$\hfill\square$

\vsp

\begin{rem} As constructed above, the dynamical realization depends
not only on the left-ordering $\preceq$, but also on the enumeration
$(g_i)_{i\geq 0}$. Nevertheless, it is not hard to check that
dynamical realizations associated to different enumerations (but the
same ordering) are \textit{topologically conjugate}.\footnote{Two
actions $\phi_1\!: \Gamma \to \mathrm{Homeo}_+(\R)$ and
$\phi_2\!:\Gamma \to \mathrm{Homeo}_+(\R)$ are topologically
conjugate if there exists $\varphi \in \mathrm{Homeo}_+(\R)$ such
that $\varphi\circ \phi_1(g) = \phi_2(g) \circ \varphi$ for all $g
\in \Gamma$.} Thus, up to topological conjugacy, the dynamical
realization depends only on the ordering $\preceq$ of $\Gamma$.

\vsp

An important property of dynamical realizations is that they do not
admit global fixed points (\textit{i.e.,} no point is stabilized by
the whole group). Another important property is that $g\succ id$ if
and only if $g(t(id))> t(id)$, which allows us to recover the
left-ordering from the dynamical realization.
\end{rem}

\vsp

\noindent \textbf{Proof of Theorem \ref{laprop}: } First fix $a\in
G_1$ and $b\in G$ exactly as above, that is, such that
$bab^{-1}=a^r$, where $r\in \Q$, $r>1$, and $\varphi(a)(x)=x+1$,
$\varphi(b)(x)=rx$. Now let $\preceq$ be a left-ordering on $G$, and
consider its dynamical realization.  To prove Theorem \ref{laprop},
we will distinguish two cases:

\vs

\noindent \textbf{Case 1.} The element $a \in G$ is cofinal (that
is, for every $g \in G$, there are $n_1,\, n_2 \in \mathbb{Z}$ such
that $a^{n_1}\prec g\prec a^{n_2}$).

\vs

Note that in a Conradian ordering $G_1$ is convex. So $a$ can not be
cofinal. Thus, in this case we have to prove that $\preceq$ is an
Smirnov-type ordering.

\vsp

For the next two claims, recall that for any measure $\mu$ on a
measurable space $X$ and any measurable function $f: X\to X$, the
{\em push-forward measure} $f_*(\mu)$ is defined by $f_*(\mu)(A) =
\mu(f^{-1}(A))$, where $A\subseteq X$ is a measurable subset. Note
that $f_*(\mu)$ is trivial if and only if $\mu$ is trivial.
Moreover, one has $(fg)_*(\mu) = f_*(g_*(\mu))$ for all measurable
functions $f,g$.

\vsp

Similarly, the {\em push-backward measure} $f^*(\mu)$ is defined by
$f^*(\mu)(A)=\mu(f(A))$.

\vsp\vsp\vsp\vsp

\noindent{\underbar{Claim 1.}} The subgroup $G_1$ preserves a Radon
measure $\nu$ ({\em i.e.,} a measure which is finite on compact
sets)  on the real line which is unique up to scalar multiplication
and has no atoms.

\vsp\vsp\vsp

Since $a$ is cofinal and $G_1$ is rank-one Abelian, its action on
the real line is {\em free} (that is, no point is fixed by any
non-trivial element of $G_1$). By H\"older's theorem (see
\cite[Theorem 6.10]{ghys} or \cite[\S 2.2]{book}), the action of
$G_1$ is semi-conjugated to a group of translations. More precisely,
there exists a non-decreasing, continuous, surjective function $\rho
\!: \mathbb{R} \rightarrow \mathbb{R}$ such that, to each $g \in
G_1$, one may associate a translation parameter $c_g$ so that, for
all $x \in \mathbb{R}$,
$$ \rho(g(x))=\rho(x)+ c_g.$$
Now since the Lebesgue measure $Leb$ on the real line is invariant
by translations, the {\em push-backward measure} $ \nu =
\rho^*(Leb)$ is invariant by $G_1$. Since $Leb$ is a Radon measure
without atoms, this is also the case for $\nu$.

\vsp

To see the uniqueness of $\nu$ up to scalar multiple we follow
\cite[\S 2.2.5]{book}. Given any measure $\mu$, invariant by the
action (in this case) of $G_1$, we define the associated {\em
translation number homomorphism} $\tau_\mu:G_1 \to \R$, by
$$ \tau_{\mu}(g)= \left\{ \begin{array}{c c} \mu([x,g(x)]) &
\text{ if $g(x) > x$, }  \\0 & \text{if } g(x)=x, \\    -\mu([g(x),
x]) & \text{ if $g(x) < x. $}
\end{array} \right.$$
One easily checks that this definition is independent of $x\in \R$,
and that the kernel of $\tau_\mu$ coincides with the elements having
fixed points, which in this case is just the identity of $G_1$. Now,
by \cite[Proposition 2.2.38]{book}, to prove the uniqueness of
$\nu$, it is enough to show that, for any non-trivial $\mu$,
$\tau_\mu(G_1)$ is dense in $\R$. But since $G_1$ is rank-one
Abelian, and $G_1\not\simeq \Z$, any non-trivial homomorphism from
$G_1$ to $\R$ has a dense image. In particular $\tau_\mu(G_1)$ is
dense in $\R$. So Claim 1 follows.

\vs

\noindent{\underbar{Claim 2.}} For some $\lambda\not= 1$, we have
$b_*(\nu) = \lambda \nu$.

\vsp\vsp\vsp

Since $G_1 \lhd G$, for any $a' \in G_1$ and all measurable $A
\subset \mathbb{R}$ we must have
$$b_*(\nu)(a' (A)) = \nu(b^{-1}a'
(A))=\nu(\bar{a}(b^{-1}(A)))=\nu(b^{-1}(A))=b_*(\nu)((A))$$ for some
$\bar{a} \in G_1$. (Actually, $a' = \bar{a}^r$.) Thus $b_*(\nu)$ is
a measure that is invariant by $G_1$. The uniqueness of the
$G_1$-invariant measure up to scalar factor yields $b_*(\nu)=\lambda
\nu$ for some $\lambda > 0$. Assume for a contradiction that
$\lambda$ equals 1. Then the whole group $G$ preserves $\nu$. In
this case, there is a {\em translation number homomorphism}
$\tau_{\nu} \!: G\to \mathbb{R}$ defined by
$$ \tau_{\nu}(g)= \left\{ \begin{array}{c c} \nu([x,g(x)]) &
\text{ if $g(x) < x$, }  \\ 0 & \text{if } g(x)=x, \\
-\nu([g(x), x]) & \text{ if $g(x) < x. $}
\end{array} \right .$$
The kernel of $\tau_{\nu}$ must contain the commutator subgroup of
$G$, and, since $a^{r-1}=[a,b] \in [G,G]$, we have that
$\tau_{\nu}(a^{r-1})=0$, hence $\tau_{\nu}(a)=0$. Nevertheless, this
is impossible, since the kernel of $\tau_{\nu}$ coincides with the
set of elements having fixed points on the real line (see \cite[\S
2.2.5]{book}). So Claim 2 is proved.

\vsp\vsp\vsp

By Claims 1 and 2, for each $g \in G$ we have $g_*(\nu) = \lambda_g
(\nu)$ for some $\lambda_g > 0$. Moreover, $\lambda_a = 1$ and
$\lambda_b = \lambda\not=1$. Note that, as
$(fg)_*(\nu)=f_*(g_*(\nu))$, the correspondence $g\to \lambda_g$ is
a group homomorphism from $G$ to $\R_+$, the group of positive real
numbers under multiplication. Since $G_1$ is in the kernel of this
homomorphism and any $g\in G$ is of the from $b^na^s$ for $n\in \Z,
\; s\in \Q$, we have that the kernel of this homomorphism is exactly
$G_1$.

\vsp

\begin{lem}\label{lema A} {\em Let $A:G\to Af_+(\R)$, $g\to A_g$, be defined by
$$\label{afin}A_g(x)=\esp \frac{1}{\lambda_g} x +
\frac{sgn(g)}{\lambda_g}\, \nu([t(g^{-1}),t(id)]),$$ where
$sgn(g)=\pm 1$ is the sign of $g$ in $\preceq$ (that is, $sgn(g)=1$ if $g$ is non-negative, and $sgn(g)=-1$ if $g$ is negative.). Then $A$ is an
injective homomorphism. }
\end{lem}

\noindent {\em Proof:} For $ g,h\in G$ both positive in $\preceq$,
we compute
\begin{eqnarray*}
A_{gh}(x)
 &=& \frac{1}{\lambda_{gh}} x +
\frac{1}{\lambda_{gh}} \nu([t((gh)^{-1}),t(id)]) \\
 &=&   \frac{1}{\lambda_{g}\lambda_h} x +\frac{1}{\lambda_{g}\lambda_h}
 \left[ (h_*\nu) ([t(g^{-1}),t(h)]) \right]   \\
 &=&   \frac{1}{\lambda_{g}\lambda_h} x +\frac{1}{\lambda_{g}\lambda_h}
 \left[ \lambda_h\nu ([t(g^{-1}),t(id)]) +\nu ([t(h^{-1}),t(id)])  \right]   \\
 &=&  \frac{1}{\lambda_{g}\lambda_h} x +\frac{1}{\lambda_{g}}\nu
 ([t(g^{-1}),t(id)])+
 \frac{1}{\lambda_g\lambda_h}\nu ([t(h^{-1}),t(id)]) \\
 &=& A_g(A_h(x)).
 \end{eqnarray*}

\vsp

The other cases can be treated analogously.

\vsp

Now, assume that $A_g(x)=x$ for some non-trivial $g\in G$. Then
$\lambda_g=1$. In particular $g\in G_1$, since the kernel of the
application $g\to \lambda_g$ is $G_1$. But in this case we have that
$g$ has no fixed point, so assuming that
$0=\lambda_g^{n-1}\nu([t(g^{-1}),t(id)]=\nu([t(g^{-n}),t(id)]$
implies $\nu$ is the trivial measure. This contradiction settles
Lemma \ref{lema A} . $\hfill\square$

\vs

Now, for $x\in \R$, let $F(x) = sgn(x - t(id)) \cdot
\nu([t(id),x])$. (Note that $F(t(id)) = 0$.) By semi-conjugating the
dynamical realization by $F$ we (re)obtain the faithful
representation $A \!: G \to Af_+(\R)$. More precisely, for all $g
\in G$ and all $x \in \mathbb{R}$ we have \begin{equation}
\label{semiconj} F (g(x)) = A_g (F(x))\end{equation}

For instance, if $x > t(id)$ and $g \succ id$, then
\begin{eqnarray*}
F (g(x))
 &=& \nu ([t(id),g(x)])\\
 &=& \frac{1}{\lambda_g} \nu ([t(g^{-1}),x])\\
 &=& \frac{1}{\lambda_g} \nu([t(g^{-1}),t(id)]) + \frac{1}{\lambda_g} \nu([t(id),x])\\
&=& \frac{1}{\lambda_g} F(x) + \frac{1}{\lambda_g}
\nu([t(g^{-1}),t(id)]).
\end{eqnarray*}

The action $A$ induces a (perhaps partial) left-ordering
$\preceq_A$, namely $g \succ_A id$ if and only if $A_g (0) > 0$.
Note that equation (\ref{semiconj}) implies that for every $g\in
G_1$, $g\succ id$, we have $A_g(0)>0 $ so $g\succ_A id$, and for
every $f\in G$ such that $A_f(0)>0$, we have $f\succ id$. In
particular, if the orbit under $A$ of $0$ is free (that is, for
every non-trivial element $g\in G$, we have $A_g(0)\not=0$), then
(\ref{semiconj}) yields that $\preceq_A$ is total and coincides with
$\preceq$ (our original ordering).

\vsp

If the orbit of $0$ is not free (this may arise for example when $\,
t(id) \,$ does not belong to the support of $\nu$), then the
stabilizer of $0$ under the action of $A$ is isomorphic to $\Z$.
Therefore, $\preceq$ coincides with either $\preceq_A^+$ or
$\preceq_A^-$ (the definition of $\preceq_A^{\pm}$ is similar to
that of $\preceq_{\varepsilon}^{\pm}$ above).

\vsp

At this point we have that $\preceq$ can be realized as an induced
ordering from the action given by $A$. Therefore arguing as in the
proof of Proposition \ref{C y S} we have that $\preceq_A$, and so
$\preceq$, is non-isolated.

\vsp

To show that $\preceq$ is an Smirnov-type ordering, we need to
determine all possible embeddings of $G$ into the affine group.
Recall that $bab^{-1}=a^r$, $r=p/q>1$.

\vsp

\begin{lem} {\em Every faithful representation of $\,G$
in the affine group is given by}
$$ a\sim \left( \begin{array}{c c} 1 & \alpha \\ 0 &1 \end{array}
\right), \;\;\; b\sim \left( \begin{array}{c c} r&\beta  \\ 0&1
\end{array} \right)$$
for some $\alpha \not=0$ and $\beta \in \mathbb{R}$.
\end{lem}

\noindent \textit{Proof:} Arguing as in Lemma \ref{lema afin} one
may check that $\varphi^\prime_{a,b}:\{a,b\}\to Af_+(\R)$ defined by
$\varphi_{a,b}^\prime(a)(x)=x+\alpha$ and
$\varphi^\prime_{a,b}(b)(x)=rx+\beta$ may be (uniquely) extended to
an homomorphic embedding $\varphi_{a,b}^\prime:G\to Af_+(\R)$.
Conversely, let
$$a\sim \left( \begin{array}{c c} s & \alpha \\ 0 &1 \end{array}
\right),\;\;\; b\sim\left( \begin{array}{c c} t & \beta \\ 0 &1
\end{array} \right)$$
be a representation. Since we are dealing with orientation
preserving affine maps, $s,t$ are positive real numbers. Moreover,
the following equality must hold:
$$a^p\sim \left( \begin{array}{c c} s^p & s^{p-1}\alpha +\ldots +s \alpha+\alpha \\ 0 &1 \end{array}
\right)= \left( \begin{array}{c c} s^q & s^{q-1}\alpha t+ s^{q-2}\alpha t +\ldots+ \alpha t - s^q\beta+\beta \\
0 & 1 \end{array} \right)\sim ba^qb^{-1}.$$ Thus $s = 1$, $t =
p/q=r$. Finally, since the representation is faithful,
$\alpha\not=0$. $\hfill \square$

\vspace{0.35cm}

Let $\alpha,\beta$ be such that $A_a(x)=x+\alpha$ and
$A_b(x)=rx+\beta$. We claim that if the stabilizer of $0$ under $A$
is trivial --which implies in particular that $\beta \!\neq\! 0$-- ,
then $\preceq_A$ (and hence $\preceq$) coincides with
$\preceq_{\varepsilon}$ if $\alpha > 0$ (resp.
$\overline{\preceq}_{\varepsilon}$ if $\alpha < 0$), where
$\varepsilon = \frac{\beta}{(r-1) \alpha} $. Indeed, if $\alpha >
0$, then for each $g = b^n a^s \in G$, $s\in \Q$, we have $A_g(0) =
r^n s\alpha +\beta \frac{r^n - 1}{r-1} $. Hence $A_g (0) > 0$ holds
if and only if
$$r^n \frac{\beta}{(r-1)\alpha } +r^n s >\frac{\beta}{(r-1)\alpha}.$$
Letting $\varepsilon := \frac{\beta}{(r-1)\alpha}$, one easily
checks that the preceding inequality is equivalent to
$g\succ_\varepsilon id$. The claim now follows.

\vsp

In the case where the stabilizer of $0$ under $A$ is isomorphic to
$\mathbb{Z}$, similar arguments to those given above show that
$\,\preceq \,$ coincides with either $\,\preceq_{\varepsilon}^+$, or
$\, \preceq_{\varepsilon}^{-}$, or $\,
\overline{\preceq}_{\varepsilon}^+ \,$, or $\,
\overline{\preceq}_{\varepsilon}^{-} \,$, where $\varepsilon$ again
equals $ \frac{\beta}{(r-1) \alpha}$.

\vspace{0.43cm}

\noindent\textbf{Case 2.} The element $a\in G$ is not cofinal.

\vsp

In this case, for the dynamical realization of $\,\preceq \,, \,$
the set of fixed points of $a$, denoted $Fix(a)$, is non-empty. We
claim that $b(Fix(a))=Fix(a)$. Indeed, let $r=p/q$, and let $x\in
Fix(a)$. We have
$$a^{p}(b(x))=a^{p}b(x)=ba^q(x)=b(x)\, .$$
Hence $a^p(b(x))=b(x)$, which implies that $a(b(x)) = b(x)$ as
asserted. Observe that since there is no global fixed point for the
dynamical realization, we must have $b(x)\not=x \,, \,$ for all $x
\in Fix(a) \,.$ Note also that, since $G_1$ is rank-one Abelian
group, $Fix(a)=Fix(G_1)$.

\vsp

Now let $x_{-1} = \inf \{t(g)  \mid g\in G_1\}$ and $x_1=\sup
\{t(g)\mid g\in G_1\}$. It is easy to see that $x_{-1}$ and $x_1$
are fixed points of $G_1$. Moreover, $x_{-1}$ (resp. $x_1$) is the
first fixed point of $a$ on the left (resp. right) of $t(id)$. In
particular, $b((x_{-1},x_1))\cap(x_{-1},x_1)=\emptyset$, since
otherwise one may create a fixed point inside $(x_{-1},x_1)$. Taking
the {\em reverse} ordering if necessary, we may assume $b\succ id$.
In particular, we have that $b(x_{-1})\geq x_1$.

\vsp

We now claim that $G_1$ is a convex subgroup. First note that, by
the definition of the dynamical realization, for every $g\in G$ we
have $\,t(g)=g(t(id))$. Then, it follows that for every $g\in G_1$,
$t(g) \!\in (x_{-1},x_1)$. Now let $m,s \in \mathbb{Z}$ and $g\in
G_1$ be such that $id \prec b^m g \prec a^s$. Then we have $t(id) <
b^m (t(g)) < t(a^s) < x_1$. Since $b(x_{-1}) \geq x_1$, this easily
yields $m = 0$, that is, $b^m g =g\in G_1$.

\vsp

We have thus proved that $G_1$ is a convex (normal) subgroup of $G$.
Since the quotient $G / G_1$ is isomorphic to $\mathbb{Z}$, an
almost direct application of Theorem \ref{teo C} shows that the
ordering $\, \preceq \,$ is Conradian. This concludes the proof of
Theorem \ref{laprop}. $\hfill\square$

\vsp

\begin{rem} It follows from Theorem \ref{laprop} and Proposition \ref{C y S}
that no left-ordering is isolated in $\mathcal{LO}(G)$. Therefore,
since any group with normal rational series is countable,
$\mathcal{LO}(G)$ is a totally disconnected Hausdorff and compact
metric space, thus homeomorphic to the Cantor set.

%Moreover, the natural conjugacy action of $G$ on $\mathcal{LO}(G)$
%is `almost' transitive, see REFERENCIA FOR THE ACTTION. More
%precisely, for any irrational $\varepsilon \,, \,$ the orbit of $\,
%\preceq_\varepsilon \,$ under $G$ is dense in the subspace $V_a$
%formed by the left-orderings for which $a$ is positive. This easily
%follows from the fact that, for all $g\in G$, we have
%$g(\preceq_\varepsilon) = \esp \preceq_{g^{-1}(\varepsilon)}\!\!.$
%The complementary subspace $V_{a^{-1}}$ of $\mathcal{LO}(G)$ is
%densely covered by the orbit of the reverse ordering of
%$\;{\preceq}_\varepsilon.$
\end{rem}

\begin{rem} The above method of proof also gives a complete classification --up to
topological semiconjugacy-- of all actions of $G$ by
orientation-preserving homeomorphisms of the real line (compare
\cite{Na-sol}). In particular, all these actions come from
left-orderings on the group (compare with Question 2.4 in
\cite{navas} and the comments before it).
\end{rem}

%%%%%%%%%%%%%%%%%%%%%%%%%%%%%%%%%%%%%%%%%%%%%%%%%%%%%%%

\section{The general case}
\subsection{A technical proposition}

%%%%%%%%%%%%%%%%%%%%%%%%%%%%%%%%%%%%%%%%%%%%%%%%%%%%%%%

\hspace{0.4 cm} The main objective of this section is to prove the
following

\begin{prop} \label{max}\textit{ Let $G$ be a group with only finitely many $C$-orderings,
and let $H$ be its maximal convex subgroup (with respect to any
$C$-ordering). Then $H$ is a Tararin group, that is, a group with
only finitely many left-orderings.}
\end{prop}

Note that the existence of a maximal convex subgroup follows from
Theorem \ref{teo B}. Note also that Proposition \ref{max} implies
that no group with only finitely many $C$-orderings, whose rational
series has length at least 3, is bi-orderable (see also
\cite[Proposition 3.2]{rivas}).

\vsp

The proof of Proposition \ref{max} is a direct consequence of the
following

\begin{lem} \label{useful}\textit{Let $G$ be a group with only finitely many
$C$-orderings whose rational series has length at least three:
\begin{equation}\label{serie convexa} \{id\}=G_0\lhd G_1\lhd G_2\lhd \ldots \lhd G_n=G\, ,\;\; n\geq
3. \end{equation} Then given $a\in G_1$ and $b\in G_i$, $ i\leq
n-1$, we have that $bab^{-1}= a^{\varepsilon}$, $\varepsilon=\pm
1$.}
\end{lem}

\noindent \textit{Proof:} We shall proceed by induction on $i$. For
$i=0,1$ the conclusion is obvious. We work the case $i=2$. Let $b\in
G_2$, and suppose that $bab^{-1}=a^r$, where $r\not=\pm 1$ is
rational. Clearly this implies that $b^nab^{-n}=a^{r^n}$ for all
$n\in \Z$.

\vsp

Since $G_3/G_1$ is non-Abelian, there exists $c\in G_3$ such that
$cb^pc^{-1}=b^q w$, with $p\not=q$ integers and $w\in G_1$. Note
that $wa=aw$. We let $t\in \Q$ be such that $cac^{-1}=a^t$. Then we
have
$$a^{r^q}=b^qab^{-q}=b^q\,waw^{-1} b^{-q}=cb^pc^{-1} a\,
cb^{-p}c^{-1}= cb^p a^{1/t} b^{-p}c^{-1}= c
a^{\frac{r^p}{t}}c^{-1}=a^{r^p},$$ which is impossible since $r\not=
\pm 1$ and $p\not= q$. Thus the case $i=2$ is settled.

Now assume, as induction hypothesis, that for any $w\in G_{i-1}$ we
have that $waw^{-1}=a^\varepsilon$, $\varepsilon=\pm 1$. Suppose
also that there exists $b\in G_i$ such that $bab^{-1}=a^r$,
$r\not=\pm1$. As before, we have that $b^nab^{-n}=a^{r^n}$ for all
$n\in \Z$.

\vsp

Let $c\in G_{i+1}$ such that $cb^pc^{-1}=b^q w$, with $p\not=q$
integers and $w\in G_{i-1}$. Let $t\in \Q$ be such that
$cac^{-1}=a^t$. Then we have
$$a^{r^q}= b^qab^{-q}=b^q\,w\, w^{-1}aw\, w^{-1} b^{-q}=cb^pc^{-1} a^{\varepsilon}\,
cb^{-p}c^{-1}= cb^p a^{\varepsilon/t} b^{-p}c^{-1}= c
a^{\frac{\varepsilon r^p}{t}}c^{-1}=a^{\varepsilon r^p},$$ which is
impossible since $r\not= \pm 1$ and $p\not= q$ implies
$|r^p|\not=|r^q|$. This finishes the proof of Lemma \ref{useful}.
$\hfill\square$

\vs

\noindent {\em Proof of Proposition \ref{max}:} Since in any
Conradian ordering of $G$, the convex series is precisely the
rational series, we have that $H=G_{n-1}$ in (\ref{serie convexa}).
So $H$ has a rational normal series. Therefore, to prove that $H$ is
a Tararin group, we only need to check that no quotient
$G_{i}/G_{i-2}$, $2\leq i\leq n-1$, is bi-orderable.

\vsp

Now, if in (\ref{serie convexa}) we take the quotient by the normal
and convex subgroup $G_{i-2}$, Lemma \ref{useful} implies that
certain element in $G_{i-1}/G_{i-2}$ is sent into its inverse by the
action of some element in $G_{i}/G_{i-2}$. Thus $G_{i}/G_{i-2}$ is
non-bi-orderable.$\hfill\square$

\begin{cor} \label{buen remark} \textit{A group $G$ having only finitely many
$C$-orderings, with rational series
$$\{id\}\lhd G_1\lhd\ldots\lhd G_{n-1}\lhd G_n=G,$$
is a Tararin group if and only if $G/G_{n-2}$ is a Tararin group.}
\end{cor}

\vsp

%%%%%%%%%%%%%%%%%%%%%%%%%%%%%%%%%%%%%%%%%%%%%%%%%%%%%%%%%%%%%%

\subsection{Proof of the Main Theorem}

%%%%%%%%%%%%%%%%%%%%%%%%%%%%%%%%%%%%%%%%%%%%%%%%%%%%%%%%%%%%%%%

\hspace{0.4 cm} Let $G$ be a group with rational series
$$\{ id \} = G_0 \lhd G_{1} \lhd \ldots \lhd G_{n-1} \lhd G_n = G,\;\; n\geq3,$$
such that no quotient $G_i/G_{i-2}$ is Abelian. Moreover, assume $G$
is not a Tararin group. Let $\preceq$ be a left-ordering on $G$. To
show that $\preceq$ is non-isolated we will proceed by induction.
Therefore, we assume as induction hypothesis that no group with only
finitely many $C$-orderings, but infinitely many left-orderings,
whose rational series has length less than $n$, has isolated
left-orderings.

\vsp

The main idea of the proof is to find a convex subgroup $H$ such
that either $H$ has no isolated left-orderings or such that $H$ is
normal and $G/H$ has no isolated left-orderings. Indeed, by Corollary \ref{extension}, this is enough to show that $\preceq$ is non-isolated. We will see that the appropriate convex subgroup to look at is the {\em convex closure of $G_1$} (with respect to $\preceq$), that is, the smallest
convex subgroup that contains $G_1$.

\vsp

For $x,y\in G$, consider the relation in $G$ given by $x\sim y $ if
and only if there are $g_1,g_2\in G_1$ such that $g_1x\preceq y
\preceq g_2x$. We check that $\sim$ is an equivalence relation.
Clearly $x\sim x$ for all $x\in G$. If $x\sim y$ and $y\sim z$ then
there are $g_1,g_2,g_1^\prime, g_2^\prime\in G_1$ such that
$g_1x\preceq y \preceq g_2x$ and $g_1^\prime y \preceq z\preceq
g_2^\prime y$. Then $g_1^\prime g_1x\preceq z\preceq g_2^\prime g_2
x$, so $x\sim z$. Finally $g_1x\preceq y \preceq g_2x$ implies
$g_2^{-1} y \preceq x\preceq g_1^{-1} y$, so $x\sim y $ implies $
y\sim x$.

\vsp

Now let $g,x,y$ in $G$ be such that $x\sim y$, hence $g_1x\preceq y
\preceq g_2x$, for some $g_1,g_2 \in G_1$. Then $gg_1x\preceq
gy\preceq gg_2x$. Since $G_1$ is normal we have that
$gg_1x=g_1^\prime gx$ and $ gg_2x=g_2^\prime g x$, for some
$g_1^\prime, g_2^\prime \in G_1$. Therefore, $g_1^\prime gx\preceq
gy \preceq g_2^\prime gx$, so $ gx\sim gy$. That is, $G$ preserves
the equivalence relation $\sim$. Let $H:=\{x\in G\mid x\sim id\}$.

\vs

\noindent \underbar{Claim 1:} For every $g\in G$ we have
$$gH\cap H=\left\{ \begin{array}{c l} \emptyset & \text{if $g\notin H$,}
\\ H& \text{if $g\in H$}.  \end{array}
\right. $$

\vsp

Indeed, if $g\in H$, then $g\in (gH\cap H)$. Now, since $x\sim
y$ if and only if $ gx\sim gy$, we have that $gH=H$. Now suppose $g$ is
such that there is $z\in gH\cap H$. Then $id \sim z\sim g$, which
implies $g\in H$. So Claim 1 follows.

\vsp

Claim 1 implies that $H$ is a convex subgroup of $G$ that contains
$G_1$. Moreover, we have

\vs

\noindent \underbar{Claim 2:} The subgroup $H$ is the convex closure
of the subgroup $G_1$.

\vsp

Indeed, let $C$ denote the convex closure of $G_1$ in $\preceq$.
Then $H$ is a convex subgroup that contains $G_1$. Thus $C\subseteq
H$.

\vsp

To show that $H\subseteq C$ we just note that, by definition, for
every $h\in H$, there are $g_1,g_2\in G_1$ such that $g_1\preceq h
\preceq g_2$. So $H\subseteq C$, and Claim 2 follows.

\vsp

Proceeding as in Lemma \ref{lema 1} we conclude that there exists $c\in
G$ such that $c\,G_{n-1}$ generates the quotient $G/G_{n-1}$. We
have

\vs

\noindent \underbar{Claim 3:} $H/G_1$ is either trivial or
isomorphic to $\Z$.

\vsp

By proposition \ref{max}, $G_{n-1}$ is a Tararin group. Therefore, in
the restriction of $\preceq$ to $G_{n-1}$, $G_1$ is convex. So we
have that $H\cap G_{n-1}=G_1$. This means that for every $g\in
G_{n-1}\setminus G_1$, $gH\cap H=\emptyset$.

\vsp

Now, assume $H/G_1$ is non-trivial and let $g\in H\setminus G_1$. By
the preceding paragraph we have that $g\notin G_{n-1}$. Therefore,
$g=c^{m_1} w_{m_1}$, for $m_1\in \Z$, $m_1\not=0$ and $w_{m_1}\in
G_{n-1}$.

Let $m_0$ be the least positive $m\in \Z$ such that $c^mw_m\in H$,
for $w_m\in G_{n-1}$. Then, by the minimality of $m_0$, we have that
$m_1$ is a multiple of $m_0$, say $km_0=m_1$. Letting
$(c^{m_0}w_{m_0})^k=c^{m_0k}\overline{w_{m_0}}$, we have that
$(c^{m_0}w_{m_0})^{-k}c^mw_m=\overline{w_{m_0}}^{-1}w_m\in H $.
Since $\overline{w_{m_0}}^{-1}w_m\in G_{n-1}$, we have that
$\overline{w_{m_0}}^{-1}w_m\in G_1$. So we conclude that $(c^m_0
w_{m_0})^k \,G_1=c^mw_m\, G_1$, which proves our Claim 3.

\vs

We are now in position to finish the proof of the Main Theorem.
According to Claim 3 above, we need to consider two cases.

\vs

\noindent \textbf{Case 1:} $H=G_1$.

\vsp

In this case, $G_1$ is a convex normal subgroup of $\preceq$ and,
since by induction hypothesis $G/G_1$ has no isolated
left-orderings, $\preceq$ is non-isolated.

\vs

\noindent \textbf{Case 2:} $H/G_1\simeq \Z$.

\vsp

In this case, $H$ has a rational series of length 2:
$$\{id\}=G_0\lhd G_1 \lhd H.$$
We let $a\in G_1$, $a\not=id$, and $h\in H$ be such that $hG_{1}$
generates $H/G_{1}$. Let $r\in \Q$ be such that $hah^{-1}=a^r$. We
have three subcases:

\vsp

\noindent \textit{Subcase i)} $r<0$.

\vsp

Clearly, in this subcase, $H$ is non-bi-orderable. So $H$ is a
Tararin group and $G_1$ is convex in $H$. But, as proved in Claim 2,
$H$ is the convex closure of $G_1$. Therefore, this subcase does not
arise.

\vsp

\noindent \textit{Subcase ii)} $r>0$.

\vsp

Since $r>0$, we have that $H$ is not a Tararin group, thus $H$ has
no isolated left-orderings. Therefore $\preceq$ is non-isolated.

\vsp

\noindent \textit{Subcase iii)} $r=0$.

\vsp

In this case, $H$ is a rank-two Abelian group, so it has no isolated
orderings. Hence $\preceq$ is non-isolated. This finishes the proof
of the Main Theorem.

\vsp

%%%%%%%%%%%%%%%%%%%%%%%%%%%%%%%%%%%%%%%%%%%%%%%%

\subsection{An illustrative example}

%%%%%%%%%%%%%%%%%%%%%%%%%%%%%%%%%%%%%%%%%%%%%%%%

\hspace{0.4 cm} This subsection is aimed to illustrate the different
kinds of left-orderings that may appear in a group as above. To do
this, we will consider a family of groups with eight $C$-orderings.
We let $G(n)=\langle a, b , c\mid bab^{-1}=a^{-1}, cbc^{-1}=b^3,
cac^{-1}=a^n\rangle$, where $n\in \Z$. It is easy to see that $G(n)$
has a rational series of length three,
$$\{id\}\lhd G_1=\langle a\rangle \lhd G_2=\langle a, b\rangle\lhd
G(n) .$$ In particular, in a Conradian ordering, $G_1$ is convex and
normal.

\vsp

Now we note that $G(n)/G_1\simeq B(1,3)$, where $B(1,3)=\langle
\beta, \gamma \mid \gamma \beta\gamma^{-1}=\beta^3\rangle$ is a
Baumslag-Solitar group, and the isomorphism is given by
$c\to\gamma\,$, $\;\;b\to \beta \,$, $\;\; a \to id$. Now consider
the (faithful) representation $\varphi: B(1,3)\to Homeo_+(\R)$ of
$B(1,3) \simeq G(n)/G_1$ into $Homeo_+(\R)$ given by
$\varphi(\beta)(x)=x+1$ and $\varphi(\gamma)(x)=3x$. It is easy to
see that, if $x\in \R$, then $Stab_{\varphi(B(1,3))}(x)$ is either
trivial or isomorphic to $\Z$.

\vsp

In particular, $Stab_{\varphi(B(1,3))}(\frac{-3k}{2})=\langle
\gamma\beta^k\rangle$, where $k\in \Z$. Thus $\langle \gamma\beta^k
\rangle$ is convex in the induced ordering from $\frac{-3k}{2}$ (in
the representation given by $\varphi$). Now, using the isomorphism
$G(n)/G_1 \simeq B(1,3)$, we have induced an ordering on $G(n)/G_1$
with the property that $\langle cb^k\,G_1\rangle$ is convex. We
denote this left-ordering by $\preceq_2$. Now, extending $\preceq_2$
by the initial Conradian ordering on $G_1$, we have created an
ordering $\preceq$ on $G(n)$ with the property that $H(n)=\langle
a,cb^k\rangle$ is convex. Moreover, we have:

\vsp

- If $n=1$ and $k=0$, then $H(n)=\langle a,c \rangle \leq G(n)$ is
convex in $\preceq$ and $ca=ac$, as in Subcase $iii)$ above.

\vsp

- If $n\geq 2$, and $k=0$, then $H(n)=\langle a,c\rangle \leq G(n) $
is convex in $\preceq$ and $cac^{-1}=a^2$, as in Subcase $ii)$
above.

\vsp

- If $n\leq -1$ and $k$ is odd, then $H(n)=\langle a, cb^k\rangle
\leq G(n)$ is convex and $cb^k \, a\, b^{-k}c^{-1}= a^{-n}$ (again)
as in Subcase $ii)$ above.

\newpage

%%%%%%%%%%%%%%%%%%%%%%%%%%%%%%%%%%%%%%%%%%%%%%%%%%%%%%%%%%%%%%%%%%%%%%%%%%%%%%%%%%%%%%%%%%%%%%%%%%%%%%%%%
%%%%%%%%%%%%%%%%%%%%%%%%%%%%%%%%%%%%%%%%%%%%%%%%%%%%%%%%%%%%%%%%%%%%%%%%%%%%%%%%%%%%%%%%%%%%%%%%%%%%%%%%%

\begin{small}

\vspace{0.1cm}

%%%%%%%%%%%%%%%%%%%%%%%%%%%%%%%%%%%%%%%%%%%%%%%%%%%%%%%%%%%%%%%%%%%%%%%%%%%%%%%%%%%%%%%%%

\vspace{0.37cm}

\noindent Crist\'obal Rivas

\vsp

\noindent Dep. de Matem\'aticas, Fac. de Ciencias, Univ. de Chile

\vsp

\noindent Las Palmeras 3425, \~Nu\~noa, Santiago, Chile

\vsp

\noindent Email: cristobalrivas@u.uchile.cl

\end{small}

%%%%%%%%%%%%%%%%%%%%%%%%%%%%%%%%%%%%%%%%%%%%%%%%%%%%%%%%%%%%%%%%%%%%%%%%%%%%%%%%%%%%%%%%%%
%%%%%%%%%%%%%%%%%%%%%%%%%%%%%%%%%%%%%%%%%%%%%%%%%%%%%%%%%%%%%%%%%%%%%%%%%%%%%%%%%%%%%%%%%%

\end{document}